\documentclass[11pt]{article}
\usepackage{amstext,amssymb,amsmath,amsbsy}

\usepackage{hyperref}
\usepackage{amscd}
\usepackage{amsfonts}
\usepackage{indentfirst}
\usepackage{verbatim}
\usepackage{amsmath}
\usepackage{amsthm}
\usepackage{enumerate}
\usepackage{graphicx}
\usepackage{color}
\usepackage[OT1]{fontenc}
\usepackage[latin1]{inputenc}
\usepackage[english]{babel}
\usepackage{amssymb}
\usepackage{subfig}

\newtheorem{Theorem}{Theorem}

\newtheorem{remark}{Remark}

\newtheorem{Assumption}{Assumption}

\numberwithin{equation}{section}

\newcommand{\R}{\mathbb{R}}
\newcommand{\Beq}{\begin{equation} \left\{ \begin{array}{rcll}}
\newcommand{\Eeq}{\end{array} \right.\end{equation}}
\newcommand{\ii}{\textrm{i}}
\newcommand{\ds}{\displaystyle}
\newcommand{\beq}{\begin{equation}}
\newcommand{\eeq}{\end{equation}}
\newcommand{\usc}{u_{\rm sc}}
\newcommand{\intk}{\int_{k_1}^k}
\newcommand{\intK}{\int_{k_1}^{k_2}}
\newcommand{\Id}{\mathrm{Id}}

\newcommand{\dint}{\int}

\newcommand{\labela}{True model}
\newcommand{\labelb}{Radon transform data $(\mathcal{R} \beta)(r, \theta)$ for the true $\beta(x_1, x_2, 0)$}
\newcommand{\labelc}{$\l_B(r, \theta,  60)$}
\newcommand{\labeld}{The
reconstructed image $\mathcal{R}^{-1}l_{B}\left( r,\theta ,60\right)$}
\newcommand{\labele}{$\l_B(r, \theta,  90)$}
\newcommand{\labelf}{The
reconstructed image $\mathcal{R}^{-1}l_{B}\left( r,\theta ,90\right)$}
\newcommand{\labelg}{Reconstructed $\tau(x, x_0) - |x - x_0|$}
\newcommand{\labelh}{The
approximate (Remark \ref{Re 4}) and reconstructed (solid) calculated phase $%
\tau \left( x,x_{0}\right) $ when $x_{0}=\left( 1,0,0\right) $, $x=\left(
x_{1},x_{2},0\right) $ and $\left\vert x\right\vert =1.$ The horizontal axis indicates the angle $\widehat{x_0, x}$ (in degrees).}
\newcommand{\labeli}{The
approximate (Remark \ref{Re 4}) and reconstructed (solid) approximate moduli $%
A\left( x,x_{0}\right) $ of the total field when $x_{0}=\left( 1,0,0\right) $%
, $x=\left( x_{1},x_{2},0\right) $ and $\left\vert x\right\vert =1.$ The horizontal axis indicates the angle $\widehat{x_0, x}$ (in degrees).}
\newcommand{\labelk}{The reconstructed image by the solution of Problem 1}

\newcommand{\width}{4cm}
\newcommand{\height}{3.3cm}

\title{Nanostructures imaging via numerical solution of a 3-d inverse
scattering problem without the phase information \footnote{T\MakeLowercase{he work of first two authors was supported by the Office
of} N\MakeLowercase{aval} R\MakeLowercase{esearch grant} N00014-15-1-2330 \MakeLowercase{and the} US A\MakeLowercase{rmy} R\MakeLowercase{esearch} L\MakeLowercase{aboratory
and} US A\MakeLowercase{rmy} R\MakeLowercase{esearch} O\MakeLowercase{ffice grant} W911NF-15-1-0233.}}
\author{Michael V. Klibanov\thanks{
Department of Mathematics and Statistics, University of North Carolina at
Charlotte, Charlotte, NC 28223, USA (\texttt{mklibanv{\@@}uncc.edu}).} \and
Loc Hoang Nguyen\thanks{
Department of Mathematics and Statistics, University of North Carolina at
Charlotte, Charlotte, NC 28223, USA (lnguye50@uncc.edu). } \and Kejia Pan\thanks{%
The corresponding author. School of Mathematics and Statistics No. 932, Central South University, Changsha 410083, China (pankejia@hotmail.com).}}

\begin{document}

\maketitle
\begin{abstract}
Inverse scattering problems without the phase information arise in imaging of nanostructures whose sizes are hundreds of nanometers as well as in imaging of biological cells. The governing equation is the 3-d generalized Helmholtz equation with the unknown coefficient, which represents the spatially distributed dielectric constant. It is assumed that only the
modulus of the complex valued wave field is measured on a frequency range. The phase is not measured. Two rigorous numerical methods are tested and their performances are compared for realistic ranges of parameters. These methods are based on two reconstruction procedures, which were recently proposed in \cite{KRB,KR2}.
\end{abstract}

35R30, 35L05, 78A46

Keywords: Phaseless inverse scattering problem,  imaging of nanostructures and biological cells, rigorous numerical methods .


\section{Introduction}

\label{sec:1}

This is the first
publication in which rigorous numerical methods are applied to solve
the phaseless inverse scattering problem of the reconstruction of a
coefficient of a PDE. We also mention the work of Ivanyshyn, Kress and Serranho 
\cite{Iv1} as well as that of Ivanishin and Kress \cite{Iv2}, in which surfaces of
obstacles are numerically reconstructed from the phaseless data. 

Inverse scattering problems without the phase information occur in imaging
of structures whose sizes are of the micron range or less.\ Recall that 1
micron ($1\mu m)=10^{-6}m.$ For example, nano structures typically have
sizes of hundreds on nanometers (nm). Either optical of X-ray radiation is
used in this imaging. Since $1nm=10^{-9}m,$ the sizes of these
nanostuctures are between $0.1\mu m$ and $1\mu m.$\ Therefore, the
wavelength $\lambda $ for this imaging should be approximately in the range $\lambda \in %
\left[ 0.05,1\right] \mu m.$ Another example is in imaging of living
biological cells. Sizes of cells are between $5\mu m$ and $100\mu m$ \cite%
{PM,Bio}.  It is well known that for the micron range of wavelengths only
the intensity of the scattered wave field can be measured, and the phase
cannot be measured \cite{Dar,Die,Khach,Pet,Ruhl}. Here, the intensity is the
square modulus of the scattered complex valued wave field. 

In the standard inverse scattering problem in the frequency domain, one is
supposed to determine a scatterer given measurements of the complex valued
wave field on the boundary of the domain of interest, see, e.g. \cite%
{Hu,Is,Li,Nov1,Nov2}. Unlike this, in the current paper we reconstruct scatterers when only the modulus of
that wave field is available while the phase is unknown.

As to the uniqueness question of phaseless inverse scattering problems, the
first uniqueness theorem was proven by Klibanov and Sacks in the 1-d case 
\cite{KS}. Also, see the work of Aktosun and Sacks \cite{AS} for a follow up result. Next, uniqueness in
the 3-d case was established by Klibanov in \cite{KSIAP,AML,AA}. However,
proofs in these references are not constructive ones. This prompted Klibanov
and Romanov to work on reconstruction procedures \cite{KR,KR1,KRB,KR2,KR3}.
We also refer to two reconstruction procedures for phaseless inverse
scattering problems which were recently developed by Novikov \cite{Nov3,Nov4}%
. Procedures of \cite{Nov3,Nov4} are quite different from ones of \cite%
{KR,KR1,KRB,KR2,KR3}.

We compare here performances of two numerical methods.\ The first one is
based on the Born approximation, as proposed in \cite%
{KRB},\ when the wave field is governed by a linearization of the generalized
Helmholtz equation at the standard one. However, this approximation fails for large values of
the frequency $k$ which are actually used both in \cite{KRB} and in this paper. Because of this, Klibanov and Romanov have proposed the
second linearization idea, which does not fail for large values of $k$ \cite%
{KR2}. Furthermore, the linearization in this case is done only on the last
step: to solve the so-called Inverse Kinematic\ Problem, whereas the phase is recovered without a linearization. Indeed, even though the Inverse Kinematic
Problem is investigated quite well by now, see, e.g. books \cite{LRV,R1,R2}, it is yet
unclear how to solve it numerically without the linearization.

 So, our second method is based on the
reconstruction procedure of \cite{KR2}. Still, we have a significantly new element here, compared with \cite{KR2}, see subsection \ref{subsection: 4.2}. While the case of the point sources
was considered in \cite{KR2}, in \cite{KR3} the case of incident plane waves
was considered. Two other reconstruction procedures, which are suitable for
a simpler case of the Schr\"{o}dinger equation, were developed by Klibanov
and Romanov in \cite{KR,KR1}.

In all five above listed reconstruction procedures of \cite%
{KR,KR1,KRB,KR2,KR3}, as well as in this paper, the last step is the well
known inverse Radon transform, see, e.g. the book of Natterer \cite{Nat}
about this transform. And in all five this transform should be applied only
in the limiting case of $k\rightarrow \infty .$ We, however, work with
specific values of the dimensionless frequency $k$
which occur in imaging of
nanostructures and biological cells, see Section \ref{sec: Problem settings} for our values of $k$. Therefore, it is not immediately clear whether some
good results can be obtained for these practically acceptable values of $k$
if using the techniques of \cite{KRB,KR2}. 

In Section \ref{sec: Problem settings}, we
pose the inverse problems we consider in this paper. In section 3 we briefly
discuss the Lippmann-Schwinger equation, which we use to generate the data
for inverse problems. In Section \ref{sec: Inverse} we explain how do we solve our inverse
problems numerically. In section \ref{Sec: Numer} we present numerical results. Section \ref{sec: summ}
 is devoted to a discussion of our results.

\section{Problem settings} \label{sec: Problem settings}

Let $\Omega \subset \mathbb{R}^{3}$ be a bounded domain. Let $%
G=\{|x|<R\}\subset \mathbb{R}^{3}$ be the ball of the radius $R$ with the
center at $\{0\}.$ We assume that $\Omega \subset G.$ Denote $S\left( 0,R\right) =\left\{ \left\vert x\right\vert =R\right\} .$ Let $c(x),x\in 
\mathbb{R}^{3}$ be a real valued function satisfying the following
conditions 
\begin{align}
c &\in C^{15}(\mathbb{R}^{3}),  \label{1.1}\\
\quad c(x) &=1+\beta (x),  \label{1.2}\\
\beta (x)\geq 0,\text{ }\beta (x)&=0\quad \text{for }\>x\in \mathbb{R}%
^{3}\setminus \Omega .  \label{1.3}
\end{align}
The smoothness requirement imposed on the function $c(x)$ is clarified in
the proof of Theorem 1 of \cite{KR2}. The non-negativity of the function $%
\beta (x)$ in (\ref{1.3}) is due to the fact that the function $%
c(x)\geq 1$ is the spatially distributed dielectric constant in the above
mentioned applications. In other words in the medium of our interest, the
dielectric constant exceeds the one of the vacuum. The function $c(x)$
generates the conformal Riemannian metric as 
\begin{equation}
d\tau =\sqrt{c(x)}\left\vert dx\right\vert ,|dx|=\sqrt{%
(dx_{1})^{2}+(dx_{2})^{2}+(dx_{3})^{2}}.  \label{1.31}
\end{equation}%
\ Here is the assumption which we use everywhere below:

\begin{Assumption} \label{assumption 1}
We assume that geodesic lines of the metric (%
\ref{1.31}) satisfy the regularity condition, i.e. for each two points $%
x,x_{0}\in \mathbb{R}^{3}$, there exists a single geodesic line $%
\Gamma \left( x,x_{0}\right)$ connecting them.
\end{Assumption}

For $x,x_{0}\in \mathbb{R}^{3},$ let $\tau (x,x_{0})$ be the solution of the
eikonal equation, 
\begin{equation}
|\nabla \tau (x,x_{0})|^{2}=c(x),\quad \tau (x,x_{0})=O\left( \left\vert
x-x_{0}\right\vert \right) ,\text{ }\>\mathrm{as}\text{ }\>x\rightarrow
x_{0}.  \label{1.9}
\end{equation}%
Here and below all derivatives are with respect to components of the vector $%
x=\left( x_{1},x_{2},x_{3}\right) .$ Let $d\sigma $ be the euclidean arc
length of the geodesic line $\Gamma \left( x,x_{0}\right) .$ Then \cite%
{R2,R3} 
\begin{equation}
\tau \left( x, x_0\right) =\dint\limits_{\Gamma \left( x,x_{0}\right) }\sqrt{%
c\left( \xi \right) }d\sigma .  \label{1.90}
\end{equation}%
Hence, $\tau (x,x_{0})$ is the travel time between points $x$ and $x_{0}$
due to the Riemannian metric (\ref{1.31}). Due to Assumption \ref{assumption 1}, $\tau
(x,x_{0})$ is a single-valued function of both points $x$ and $x_{0}$ in $%
\mathbb{R}^{3}\times \mathbb{R}^{3}$.

Consider the equation%
\begin{equation}
\Delta u+k^{2}c(x)u=0,\text{ \ }x\in \mathbb{R}^{3},\text{\ \ }  \label{100}
\end{equation}%
where the frequency $k=2\pi /\lambda ,$ where $\lambda >0$ is the
wavelength. We now comment on the ranges of parameters used in our
computations below. We point out that we work with ranges of parameters
which are \emph{realistic} to the applications to imaging of nanostructures
of several hundreds nanometers size as well as of living biological cells.
The range of our wavelengths is 
\begin{equation}
\lambda \in \left[ 0.078, 0.126\right] \mu m.  \label{101}
\end{equation}%
Thus, in our computations $\Omega $ is the sphere of the radius $R=1\mu m.$
Let us make variables to be dimensionless in (\ref{100}). To do so,
introduce a new variable $x^{\prime }=x/1\mu m$ while keeping the same
notations as above for brevity. Then elementary operations show that
equation (\ref{100}) remains the same. Recall that the frequency $k=2\pi
/\lambda .$ Hence, we obtain that the dimensionless frequency $k$, 
\begin{equation}
k\in \left[ 50,80\right] =\left[ k_{1},k_{2}\right] .  \label{1.54}
\end{equation}%
Thus, below all functions, parameters and variables are dimensionless. The
single non-realistic parameter in our computations is the assumption that
one can conduct real measurements on the sphere of such a small radius as $%
R = 1 \mu m.$ In fact, usually such measurements are conducted on the
surface of a sphere of several centimeters size \cite{Die}. This is a
delicate issue which will be considered in our future works.

To introduce the point source in equation (\ref{100}), we modify it as 
\begin{equation}
\Delta u+k^{2}c(x)u=-\delta (x-x_{0}),\quad x\in \mathbb{R}^{3},  \label{1.4}
\end{equation}%
where $x_{0}\in \mathbb{R}^{3}$ is the source position. We assume that the
function $u(x,x_{0}, k)$ satisfies the radiation condition, 
\begin{equation}
\frac{\partial u}{\partial r}+iku=o(r^{-1})\>\text{ }\mathrm{as}\text{ }%
\>r=|x-x_{0}|\rightarrow \infty .  \label{1.5}
\end{equation}%
Denote $u_{0}(x,x_{0}, k)$ the solution of the problem (\ref{1.4}), (\ref{1.5}%
) for the case $c(x)\equiv 1.$ Then $u_{0}$ is the incident spherical wave,%
\begin{equation*}
u_{0}(x,x_{0},k)=\frac{\exp \left( -ik\left\vert x-x_{0}\right\vert \right) 
}{4\pi \left\vert x-x_{0}\right\vert }.
\end{equation*}%
Let $u_{sc}(x,x_{0},k)$ be the scattered wave, which is due to the presence
of scatterers, in which $c(x)\neq 1$. Then 
\begin{equation}
u_{sc}(x,x_{0},k)=u(x,x_{0},k)-u_{0}(x,x_{0},k)=u(x,x_{0},k)-\frac{\exp
\left( -ik\left\vert x-x_{0}\right\vert \right) }{4\pi \left\vert
x-x_{0}\right\vert }.  \label{1.6}
\end{equation}%
It was shown in \cite{KR2} that problem (\ref{1.4}), (\ref{1.5}) has
unique solution $u$. Moreover, $u$ is in $C^{16+\alpha }\left( \left\vert x-x_{0}\right\vert \geq \eta \right)$ for any $\eta > 0$ and $\alpha \in (0, 1).$ Here $%
C^{16+\alpha }$ is the H\"{o}lder space.


We model the propagation of the electric wave field in $\mathbb{R}^{3}$ by a
single equation (\ref{1.4}) with the radiation condition (\ref{1.5}) instead
of the full Maxwell's system. This modeling was numerically justified in the
work of Beilina \cite{BMM} in the case of the Cauchy problem for the
equation $c\left( x\right) v_{tt}=\Delta v$.\ It was established in \cite%
{KR2} that, given Assumption \ref{assumption 1}, the solution $u\left( x,x_{0},k\right) $ of
the problem (\ref{1.4}), (\ref{1.5}) is the Fourier transform of the
function $v\left( x,x_{0},t\right) $ with respect to $t$. It was
demonstrated in \cite{BMM} that the component of the electric wave field,
which is incident upon the medium, significantly dominates two other
components when propagating through the medium. Furthermore, the propagation
of that dominating component is well governed by the function $v$. This
conclusion was verified via accurate imaging using electromagnetic
experimental data in, e.g. Chapter 5 of the book of Beilina and Klibanov 
\cite{BK} and in \cite{BK1,BK2,TBKF1,TBKF2}. In those experimental data only
the single component of the scattered electric wave field was measured: the
one which was originally incident upon the medium. We believe that the confirmation on the experimental data has a
significant merit.

Along with the function $u_{sc}(x,x_{0},k)$ we consider its Born
approximation. The Born approximation assumes that the resolvent series for
the Lippmann-Schwinger integral equation (section 3) converges. Next, it
considers only the first term in this series and truncates the rest of
terms. In other words, it is assumed that \cite{KRB}%
\begin{equation}
u_{sc}(x,x_{0},k)\approx u_{B,sc}(x,x_{0},k)=k^{2}\dint\limits_{\Omega }%
\frac{\exp (-ik|x-\xi |)}{4\pi |x-\xi |}\beta \left( \xi \right) u_{0}\left(
\xi ,x_{0},k\right) d\xi .  \label{1.61}
\end{equation}%
Note that the the Born approximation breaks up for
large values of $k$.

In this paper we study the following two inverse problems:

\noindent\textbf{Problem 1}. \emph{Suppose that the following function }$f\left(
x, x_0, k\right) $\emph{\ is known}%
\begin{equation}
f\left( x,x_{0},k\right) =\left\vert u_{sc}(x,x_{0},k)\right\vert
^{2},\>\forall \left( x,x_{0}\right) \in S\left( 0,R\right) \times S\left( 0,R\right) ,\>\forall k\in \left[ k_{1},k_{2}\right] ,  \label{1.51}
\end{equation}%
\emph{for a certain interval }$\left[ k_{1},k_{2}\right] \subset \left(
0,\infty \right) .$
\textit{Suppose that the function $c\left( x\right) $ satisfies conditions
\eqref{1.1}-\eqref{1.3}. Determine the function $\beta \left( x\right).$
}

\noindent \textbf{Problem 2}. \emph{Approximate the right hand side of (\ref{1.51}) as 
}
\begin{equation*}
\left\vert u_{sc}(x,x_{0},k)\right\vert ^{2}\approx \left\vert
u_{B,sc}(x,x_{0},k)\right\vert ^{2}.
\end{equation*} \emph{Let the given function }$f$\emph{%
\ be the same as in (\ref{1.51}). Assume that }%
\begin{equation}
f\left( x,x_{0},k\right) =\left\vert u_{B,sc}(x,x_{0},k)\right\vert
^{2},\>\forall \left( x,x_{0}\right) \in S\left( 0,R\right) \times S\left( 0,R\right)    \label{1.52}
\end{equation}%
\textit{for some values of $k>0$. Suppose that the function $c\left( x\right) $ satisfies conditions
\eqref{1.1}-\eqref{1.3}. Determine the function $\beta \left( x\right).$}

\begin{remark} \label{remark 1}~
\begin{enumerate}
\item The function $f\left( x,x_{0},k\right) $ is the same in (\ref{1.51}) and (%
\ref{1.52}), and it is computed in numerical simulations of the data generation
procedure via the solution of
problem (\ref{1.4}), (\ref{1.5}). However, (\ref{1.52}) means that when
solving Problem 2, we assume that $u_{sc}(x,x_{0},k)\approx
u_{B,sc}(x,x_{0},k).$
\item Problem 2 is simpler than Problem 1. Thus, in Problem 2 we impose less
restrictive conditions on the function $\beta \left( x\right) .$ We assume
instead of (\ref{1.1}) that $\beta \in C^{1}\left( \mathbb{R}^{3}\right) $
and do not impose neither Assumption \ref{assumption 1} nor Assumption \ref{assumption 2} (section \ref{sec: Inverse}).
\end{enumerate}
\end{remark}

\section{Lippmann-Schwinger equation and the forward problem}

To generate the function $f$ in \eqref{1.51}, \eqref{1.52} for our numerical studies, we
need to solve the forward problem \eqref{1.4}, \eqref{1.5}. In this section we briefly
explain how we solve it via the solution of the Lippmann-Schwinger
integral equation. Since in our numerical experiments the function $\beta
\left( x\right) \neq 0$ only within small inclusions inside of the domain $%
\Omega ,$ computations via the Lippmann-Schwinger equation are fast. Recall
that by \eqref{1.6}
\beq
	\usc(x) = u(x) - \frac{e^{- \ii k |x - x_0|}}{4\pi |x - x_0|}, \quad x \in \R^3 \setminus \{x_0\}.
\nonumber
\eeq
It follows from \eqref{1.4}, \eqref{1.5} that $\usc$ satisfies
\Beq
	\Delta \usc + k^2 \usc &=& -k^2 \beta u &\mbox{in } \R^3,\\
	\ds\frac{\partial \usc}{\partial r} + \ii k \usc &=& O(r^{-1})&\mbox{as } r \to \infty.
	\label{3.2}
\Eeq
Using Theorem 8.3 in the book of Colton and Kress \cite{CK}, we derive from \eqref{3.2} the Lippmann-Schwinger equation
\beq
	\usc(x, x_0, k) = k^2\int_{\Omega} \frac{e^{-\ii k |x - x_0|}}{4 \pi |x - x_0|}  \beta(\xi) u(\xi, x_0, k) d\xi, \quad x \in \R^3 \setminus\{x_0\}.
	\label{usc solution}
\eeq
Therefore, the integral
equation for the total field $u\left( x,x_{0},k\right)$ is
\beq
	u(x, x_0, k) = \frac{e^{-\ii k |x - x_0|}}{4 \pi |x - x_0|} + k^2\int_{\Omega} \frac{e^{-\ii k |x - \xi|}}{4 \pi |x - \xi|}  \beta(\xi) u(\xi, x_0, k) d\xi, \quad x \in \R^3 \setminus\{x_0\}.
	\label{LS eqn}
\eeq
Theorems 8.3 and 8.7 of \cite{CK} guarantee that there exists unique
solution $u\in C^{2}\left( \mathbb{R}^{3}\diagdown \left\{ x_{0}\right\}
\right) $ of the integral equation \eqref{LS eqn}.\ Furthermore, this function $u$ is
the solution of problem (\ref{1.4}), (\ref{1.5}) and, therefore, it has
a higher smoothness (section \ref{sec: Problem settings}). In our computations of the forward problem
for data generation, we solve equation \eqref{LS eqn} for a discrete set of
frequencies $k\in \left[ k_{0},k_{1}\right] $, for $x\in \Omega $ and for $%
x_{0}\in S\left( 0,R\right).$ As soon as the function $u\left(
x,x_{0},k\right) $ is found for $x\in \Omega ,$ its values for point $x\in 
\mathbb{R}^{3}\diagdown \Omega $ can be easily found by substituting these
points $x$ in the right hand side of \eqref{LS eqn}. This solution produces the function $%
\left\vert u_{sc}\left( x,x_{0},k\right) \right\vert ^{2}=f\left(
x,x_{0},k\right) $ which is the data for the above Problems 1 and 2.

\section{Inverse Problems} \label{sec: Inverse}

\subsection{Two main theorems}

Our reconstruction methods are based on Theorems \ref{thm 1} and \ref{thm 2}. Following \cite%
{KR2}, let $\zeta =(\zeta _{1},\zeta _{2},\zeta _{3})$, $\zeta =\zeta
(x,x_{0})$ be geodesic coordinates of a variable point $x$ with respect to a
fixed point $x_{0}$ in the above Riemannian metric (\ref{1.31}). By 
Assumption \ref{assumption 1}, there exists a one-to-one correspondence $x\Leftrightarrow
\zeta $ for any fixed point $x_{0}$. Therefore, for any fixed point $x_{0}$
the function $\zeta =\zeta (x,x_{0})$ has the inverse function $x=g\left(
\zeta ,x_{0}\right) .$ The function $g$ determines the geodesic line \begin{equation*}
\Gamma(x,x_{0})=\{\xi :\xi =g\left( s\zeta _{0},x_{0}\right) ,s\in \lbrack 0,\tau
(x,x_{0})]\},
\end{equation*}
 where $\zeta _{0}$ is the vector which is tangent to $\Gamma
(x,x_{0})$ at the point $x_{0}$, is directed towards the point $x$ and also $%
|\zeta _{0}|=c^{-1/2}(x_{0})$.

 Consider the Jacobian $J(x,x_{0})=\det \left( \partial \zeta
/\partial x\right) .$ The mapping $x\Leftrightarrow \zeta $ is one-to-one
and continuously differentiable. Still, this does not imply that $%
J(x,x_{0})\neq 0$ for all $x,x_{0}.$ It was established in \cite{KR2} that $%
J(x,x)=1.$ Thus, we use everywhere below the following assumption:

\begin{Assumption}\label{assumption 2}
$J(x,x_{0})>0$ for all $x,x_{0}\in \mathbb{R}^{3}.$
\end{Assumption}

Define now the function $A(x,x_{0})>0$ as%
\begin{equation}
A(x,x_{0})=\frac{c(x_{0})\sqrt{J(x,x_{0})}}{4\pi \sqrt{c(x)}\tau (x,x_{0})}.
\label{1.11}
\end{equation}

\begin{Theorem}[\cite{KR2}] \label{thm 1} {Assume that conditions (\ref{1.1})-(%
\ref{1.3}) and Assumptions \ref{assumption 1}, \ref{assumption 2} hold. Fix a point }$x_{0}\in \overline{G}.$%
{\ Then for any point }$x\in \overline{G}\diagdown \left\{
x_{0}\right\} ${\ the following asymptotic behavior holds}%
\begin{equation}
u\left( x,x_{0},k\right) =A(x,x_{0})e^{-ik\tau (x,x_{0})}+O\left( \frac{1}{k}%
\right) ,k\rightarrow \infty ,  \label{1.110}
\end{equation}%
where functions $%
\tau \left( x,x_{0}\right) $ and $A\left( x,x_{0}\right) $ are given in
\eqref{1.90} and \eqref{1.11} respectively.
\end{Theorem}

Combining this theorem with (\ref{1.6}), we obtain the following asymptotic
behavior of the function $\left\vert u_{sc}\left( x,x_{0},k\right)
\right\vert ^{2}$ for $x\neq x_{0}$ as $k\rightarrow \infty $ 
\begin{multline}
\left\vert u_{sc}\left( x,x_{0},k\right) \right\vert ^{2} = A^{2}(x,x_{0}) + \frac{1}{16\pi ^{2}\left\vert x-x_{0}\right\vert ^{2}} \\-\frac{%
A(x,x_{0})}{2\pi \left\vert x-x_{0}\right\vert }\cos \left[ k\left( \tau
(x,x_{0})-\left\vert x-x_{0}\right\vert \right) \right] +O\left( \frac{1}{k}%
\right).  \label{1.12}
\end{multline}

To solve numerically Problem 2, we use Theorem \ref{thm 2}.

\begin{Theorem}[\cite{KRB}] \label{thm 2} Let conditions (\ref{1.2}), (\ref{1.3}%
) be in place and the function $\beta \in C^{1}\left( \mathbb{R}^{3}\right) 
${\ (item 2 in Remark \ref{remark 1}). Let }$u_{B,sc}(x,x_{0},k)${\ be the
function defined in (\ref{1.61}). Then the following asymptotic formula holds%
}%
\begin{equation} \label{Radon P1}
\left|u_{B,sc}(x,x_{0},k)\right|=\frac{k}{8\pi }\Bigg[\frac{1}{|x-x_{0}|}%
\int\limits_{L(x,x_{0})}\beta (\xi )d\sigma +o\Big(1\Big)\Bigg],k\rightarrow
\infty ,
\end{equation}%
{where }$L(x,x_{0})${\ is the segment of the straight line
connecting points }$x${\ and }$x_{0}.$
\end{Theorem}

In the next two subsections, we discuss numerical methods of solutions of
Problems 1, 2. We point out that since both numerical methods end up with the
inversion of the 2D Radon transform, we are interested in imaging only a 2D
cross-section of the function $\beta \left( x\right) .$ To do this, it is
not necessary to measure the function $f\left( x,x_{0},k\right) $ in (\ref%
{1.51}) and (\ref{1.52}) for all $x,x_{0}\in S\left( 0,R\right) .$ Rather,
it is sufficient to measure it for a corresponding 2D cross-section of the
sphere $S\left( 0,R\right) ,$ i.e. for $x,x_{0}\in S\left( 0,R\right) \cap
\left\{ x_{3}=a\right\} $ for a $a=const.\in \left( -R,R\right) .$ Without
further mentioning we reconstruct images below only in the central 2D
cross-section of the ball $G ,$ i.e. in the $x_{1},x_{2}$ plane$.$
We denote this cross-section as $S_{0}\left( R\right) =\left\{
x:x_{1}^{2}+x_{2}^{2}=R^{2},x_{3}=0\right\} .$ Thus, everywhere below
\textquotedblleft Radon transform" means the 2D Radon transform and $\beta $
means the function $\beta \left( x_{1},x_{2},0\right) .$

\subsection{Radon transform and Problem 2}
One
of the main tools of the numerical reconstruction of the function $\beta $
is the 2D Radon transform.
In
our computations we compute the inverse Radon transform using a tool box of
MATLAB. For every real valued function $f\in L^{2}\left( \mathbb{R}%
^{2}\right)$, the Radon transform of $g$, suggested by MATLAB, is given by
\beq
	(\mathcal{R}g) (r, \theta) = \int_{-\infty}^{\infty} g(r \cos \theta - x_2' \sin \theta, r \sin \theta + x_2' \cos \theta) d x_2',
	\nonumber
\eeq
for $ (r, \theta) \in (-\infty, \infty) \times (0, \pi).$
The function above is just the integral along the line $l_{r, \theta}.$ Here, $l_{r, \theta}$ is the line obtained by rotating the line $\{(x_1 = r, x_2): x_2 \in (-\infty, \infty)\}$ in the $x_1 x_2$ plane around the origin by the angle $\theta.$ For convenience,  we name $(r, \theta)$ as the Radon coordinate of the line $l_{r, \theta}.$ 
Denote $\mathcal{R}^{-1}$ the inverse Radon transform.

Consider now the
Born approximation. It follows from Theorem \ref{thm 2} that
\beq
	\int_{L(x, x_0)} \beta(\xi) d \sigma = \lim_{k \to \infty}l_B(r, \theta, k).
	\label{4.2}
\eeq
where
\beq
	l_B(r, \theta, k) = \frac{8 \pi}{k} |x - x_0| f(x, x_0, k)
	\label{lB}
\eeq
where $(r, \theta)$ is the Radon coordinate of the line $L(x, x_0)$ for $(x, x_0) \in S_0(R) \times S_0(R)$. The
function $l_{B}$ approximates the Radon transform of the function $\beta
\left( x_{1},x_{2},0\right) $ \ when the pair $\left( x,x_{0}\right) \in
S_{0}\left( R\right) \times S_{0}\left( R\right) .$ This leads to the
following reconstruction formula:
\beq
	\beta \approx \mathcal{R}^{-1} l_B(r, \theta, k),
	 \quad x, x_0 \in S_0(R)
\label{4.9}
\eeq where the function $f$ is defined in \eqref{1.52} and $k$ is sufficiently large. Thus,
formula \eqref{4.9} provides an approximate solution of Problem 2.  In our results displayed in Section \ref{Sec: Numer}, $k = 60, 90$.
\subsection{Problem 1} \label{subsection: 4.2}
We now
study Problem 1. We develop in this subsection a new method of the recovery
of functions $\tau \left( x,x_{0}\right) $ and $A\left( x,x_{0}\right) $ for 
$x,x_{0}\in S_0(R) $ from (\ref{1.12}). This method is
significantly different from the one proposed in \cite{KR2}. This difference
is caused by the fact that \cite{KR2} is essentially using the assumptions
that the frequency interval $k\in \left[ k_{1},k_{2}\right] $ is
sufficiently large. However, by (\ref{1.54}) our the realistic range of
dimensionless frequencies is $\left[ k_{1},k_{2}\right] =\left[ 50, 80\right]
.$ Following \cite{KR2}, we ignore the term $O\left( 1/k\right) $ in (\ref%
{1.12}). Hence, the function $f$ now is:
\beq
	f(x, x_0, k) = a(x, x_0) + b(x, x_0) \cos\left(k \alpha(x, x_0)\right), \quad k \in [k_1, k_2]
	\label{data P2}
\eeq
where 
\begin{align}
	a(x, x_0) & = A^2(x, x_0) + \frac{1}{16 \pi^2 |x - x_0|^2}\\
	b(x, x_0) &= - \frac{A(x, x_0)}{2 \pi |x - x_0|},\\
	\alpha(x, x_0) &= \tau(x, x_0) - |x - x_0|.
\end{align}

 We now propose a stable method to extract the unknown quantities $\tau$ and $A$ from $f_2$.
Fix $x, x_0 \in S(R)$. For each $k \in [k_1, k_2],$ define
\beq
	F_1(x, x_0, k) = \intk f_2(x, x_0, \kappa) d\kappa, \quad F_2(x, x_0, k) = \intk F_1(x, x_0, \kappa) d\kappa. 
\nonumber
\eeq 
For the simplicity of notations, we can temporarily ignore the dependence on $(x, x_0)$ of these functions. Assuming that $\alpha \not = 0$, we obtain
\begin{align}
	F_1(k) &= a(k - k_1) + \frac{b}{\alpha} \sin(k \alpha) - \frac{b}{\alpha} \sin(k_1 \alpha), \\
	F_2(k) &= a \frac{(k - k_1)^2}{2} - \frac{b}{\alpha^2} \cos(k \alpha) + \frac{b}{\alpha^2} \cos(k_1 \alpha) - \frac{b}{\alpha} \sin(k_1 \alpha) (k - k_1). \label{Two Anti-Derix}
\end{align}
It follows from \eqref{data P2} and \eqref{Two Anti-Derix} that
\beq
 f(k) + \alpha^2 F_2(k) = \alpha^2 a \frac{(k - k_1)^2}{2}  - \alpha b \sin(k_1 \alpha) (k - k_1) + a + b\cos(k_1 \alpha).
 \label{4.16}
\eeq
Consider 
\beq
	 \xi = (\xi_1, \xi_2, \xi_3, \xi_4) = (\alpha^2, -\alpha^2 a, \alpha b \sin(k_1 \alpha), - a - b\cos(k_1 \alpha))
	 \nonumber
 \eeq 
 as an unknown vector. We now should solve equation \eqref{4.16}
 with respect to the vector $\xi .$ We are doing this via minimizing
the following functional
 \beq
	J(\xi) = \frac{1}{2} \intK \left(\xi_1 F_2(k) + \xi_2 \frac{(k - k_1)^2}{2} + \xi_3 (k - k_1) + \xi_4 + f_2(k)\right)^2 dk.
	\nonumber
\eeq
Since equation \eqref{4.16} is linear with respect to
the vector $\xi ,$ this functional can be minimized via setting to zero
its first derivatives with respect to the components of the vector $\xi$,
assuming that the resulting matrix is not singular.
Hence, a simple calculation yields that $\xi$ solves $\mathcal{F} \xi = \mathfrak{f}$ where
\beq
	\mathcal{F} = \left(
		\begin{array}{cccc}
			\ds \intK F_2^2(k) dk & \ds \intK \frac{(k - k_1)^2}{2} F_2(k) dk & \ds \intK (k - k_1)F_2(k) dk   & \ds \intK F_2(k) dk\\
			\ds \intK \frac{(k - k_1)^2}{2} F_2(k) dk &\ds \frac{(k_2 - k_1)^5}{20} &\ds \frac{(k_2 - k_1)^4}{8} & \ds \frac{(k_2 - k_1)^3}{6}\\
			\ds \intK (k - k_1)F_2(k)dk  &\ds \frac{(k_2 - k_1)^4}{8} & \ds \frac{(k_2 - k_1)^3}{3} &\ds \frac{(k_2 - k_1)^2}{2}\\
			\ds \intK F_2(k)& \ds \frac{(k_2 - k_1)^3}{6} &\ds \frac{(k_2 - k_1)^2}{2} & k_2 - k_1
		\end{array}						
	\right) 
	\nonumber
\eeq and
\beq 
	\mathfrak{f} = - \left(
		\begin{array}{c}
			\ds \intK f_2(k)F_2(k) dk	\\
			\ds \intK \frac{(k - k_1)^2}{2}f_2(k) dk	\\
			\ds \intK (k  - k_1) f_2(k) dk	\\
			\ds \intK f_2(k) dk	
		\end{array}
	\right).
	\nonumber
\eeq

\begin{remark} 
	In our computations presented in Section \ref{Sec: Numer}, the vector $\xi$ is calculated by solving
	\beq
		(\mathcal{F}^T \mathcal{F} + \epsilon \Id) \xi = \mathcal{F}^T \mathfrak{f}
		\nonumber
	\eeq  
	where $\mathcal{F}^T$ is the transpose of $\mathcal{F}$ and $\epsilon > 0$ is a sufficiently small number.  
\end{remark}

\begin{remark}
	Recall that $\xi$ depends on $(x, x_0)$. The function 
\beq
	(x, x_0) \in S_0(R) \times S_0(R) \mapsto \sqrt{|\xi_1|} + |x - x_0|
	\nonumber
\eeq is the phase $\tau(x, x_0)$ we need to reconstruct. Having $\tau(x, x_0)$, we solve the quadratic equation \eqref{1.12} for $A(x, x_0)$. Note that this quantity might depend on $k$ and we need to take their average.
\end{remark}

To verify
the accuracy of the reconstructed functions $\tau \left( x,x_{0}\right)
,A\left( x,x_{0}\right) $ for $x,x_{0}\in S_{0}\left( R\right) ,$ we compare
them with those of the total field $u\left( x,x_{0},k\right) $ calculated
via the numerical solution of equation \eqref{LS eqn}. By \eqref{1.110} $A\left(
x,x_{0}\right) =\left\vert u\left( x,x_{0},k\right) \right\vert +O\left(
1/k\right) .$ Since the right hand side of this equality depends on $k$, we
assign an approximate value of $A\left( x,x_{0}\right) $ as%
\begin{equation}
A\left( x,x_{0}\right) \approx \frac{1}{k_{2}-k_{1}}\dint%
\limits_{k_{1}}^{k_{2}}\left\vert u\left( x,x_{0},k\right) \right\vert dk.
\label{200}
\end{equation}%
By \eqref{1.110} we have the following approximate formula at $%
k \gg 1$
\beq	
	G(x, x_0, k) = \int_{k_1}^k u(x, x_0, \kappa) d\kappa = \frac{u(x, x_0, k) - u(x, x_0, k_1)}{-\ii \tau}.
	\nonumber
\eeq
Multiplying the equations above by $\overline G(x, x_0, k)$ and integrating the resulting
equation over the interval $\left( k_{1},k_{2}\right) $ we obtain
\beq
	\tau \int_{k_1}^{k_2}|G(x, x_0, k)|^2dk = \Re\int_{k_1}^{k_2} \left(\ii (u(x, x_0, k) - u(x, x_0, k_1)) \overline G(x, x_0, k)\right)dk, \label{201}
\eeq
which yields an approximation of $\tau(x, x_0).$

\begin{remark} \label{Re 4}
Note that one might find $\tau $ by considering $\Im\left( \log
u\right) $ and again ignoring the term $O(1/k)$ in \eqref{1.110}. However, $%
\left\vert \Im\left( \log u\right) \right\vert \in \left[ n\pi ,n\pi
+2\pi \right] $ and the integer $n\geq 0$ must be chosen carefully.
Therefore, we compare in our figures our calculated values of $A\left(
x,x_{0}\right) $ and $\tau \left( x,x_{0}\right) $ with the approximate ones
given by \eqref{200} and \eqref{201}.
\end{remark}

\begin{remark}
Note that we have
approximately found both functions $\tau \left( x,x_{0}\right) $ and $%
A\left( x,x_{0}\right) $ for $x,x_{0}\in S_0(R) $ without any
linearization. However, it is unclear how to find the target function $\beta
\left( x\right) $ using these functions. Therefore, following \cite{KR2}, we
linearize the problem below.
\end{remark}

Thus, we have found the
function $\tau \left( x,x_{0}\right) $ for all $x,x_{0}\in S\left(
0,R\right) .$ The determination of the function $\beta \left( x\right) $
from the function $\tau \left( x,x_{0}\right) $ given for all $x,x_{0}\in
S_0(R) $ is called \textquotedblleft Inverse Kinematic\
Problem". As it was
stated in Introduction, it is yet unclear how to solve this problem
numerically.
Therefore, we solve below the linearized Inverse Kinematic Problem. Assume
that $||\beta ||_{C^{2}\left( \overline{\Omega }\right) }<<1.$ We now
linearize the function $\tau \left( x,x_{0}\right) $ with respect to the
function $\beta .$ This linearization can be found in Theorem 11 of Chapter
3 of \cite{LRV}, in \S 5 of Chapter 2 of \cite{R1} and in \S 4 of Chapter 3
of \cite{R2}. We obtain 
\begin{equation}
\tau \left( x,x_{0}\right) =|x-x_{0}|+\dint\limits_{L\left( x,x_{0}\right)
}\beta \left( \xi \right) d\sigma .  \label{6.7}
\end{equation}%
To be precise, one should have \textquotedblleft $\approx $" instead of
\textquotedblleft $=$" in (\ref{6.7}). Using (\ref{6.7}), we obtain
\beq
	\beta = \mathcal{R}^{-1}(\tau(x, x_0) - |x - x_0|), \quad (x, x_0) \in S_0(R) \times S_0(R). \label{reconstructed beta}
\eeq

\begin{remark}
	In the case that the line $L(x, x_0)$ passes a scatterer, we call $\usc(x, x_0, k)$ a forward scattering wave. If there is no inclusion on the line $L(x, x_0)$, $\usc(x, x_0, k)$ is called either side or back scattering wave. We have observed in our computations that forward scattering waves dominate the side and back ones. 
We believe that this is the main reason
why the inverse Radon transform, which is actually based on an analysis of
projections, provides good quality images even though we work with the
propagation of waves.	
Therefore, due to \eqref{6.7} and the reconstruction method in \eqref{reconstructed beta}, we only find $\tau$ when
	\beq
		\frac{\|\usc(x, x_0, \cdot)\|_{L^2(k_1, k_2)}}{\|\usc\|_{L^2(S_0(R) \times S_0(R) \times (k_1, k_2))}} > \epsilon
		\nonumber
	\eeq
	for some sufficiently small positive number $\epsilon$ (in our computation $\epsilon = 4\cdot 10^{-4}$).
Otherwise, we set $\tau(x, x_0) = |x - x_0|$.
This truncation enhances the resolution of the images, see item (j) in our Figures.
\end{remark}

\section{Numerical results} \label{Sec: Numer}

In the section the domain $G=\left\{ x:\left\vert
x\right\vert <1\right\} $. Hence, in this section $R=1.$ The domain $\Omega
\subset G$ is the cube inscribed in $G$ and sides of this cube are parallel
to coordinate planes. Since by \eqref{1.3} the function $\beta \left( x\right) =0$
outside of this cube $\Omega ,$ then images below are presented only in the
central 2D cross section of $\Omega ,$ i.e. in the square $\left\{
\left\vert x_{1}\right\vert ,\left\vert x_{2}\right\vert <\sqrt{2}%
/2,x_{3}=0\right\} .$ The range of dimensional frequencies $k$ is as in
\eqref{1.54}. Recall that by \eqref{101} this means that wavelengths 
$\lambda \in \left[ 0.078,0.126\right] \mu m.$ Due to the smoothness requirements imposed on the
function $\beta $ in Theorems \ref{thm 1}, \ref{thm 2}, we introduce the shape function $\varphi
\in C_{0}^{\infty }\left( \mathbb{R}^{3}\right) $ as:
\beq
	\varphi(x) = \left\{
		\begin{array}{ll}
			\ds e^{1 - \frac{1}{1 - |x|^2}} & |x| < 1, \\
			0 &\mbox{otherwise}
		\end{array}
	\right.
	\nonumber
\eeq	
to construct some true models based on $\varphi.$ 
Here
is how we construct each inclusion. First, we consider a ball of the radius $%
r$ with the center at the point $y$. And we set the initial function $\beta
_{init}\left( x\right) $ as $\beta _{init}\left( x\right) =\gamma $ inside
this inclusion and $\beta _{init}\left( x\right) =0$ outside of it. Here $%
\gamma =const.>0.$ Next, we set $\beta \left( x\right) =\beta _{init}\left(
x\right) \varphi \left( \left( x-y\right) /r\right) .$ Hence, $\gamma $ is
the maximal value of the function $\beta \left( x\right) $ in that
inclusion. If we have two inclusions, we act similarly. Below the
\textquotedblleft distance between surfaces of two inclusions" means the
distance between surfaces of those two
original spheres and the \textquotedblleft radius" of an inclusion means the
radius of that original sphere. Clearly the number $\gamma +1$ can be considered as
the inclusion/background contrast of the function $c\left( x\right) =1+\beta
\left( x\right) .$ In all our numerical examples $\gamma = 1.$

In this section, we show numerical results in several cases.  For each case we show
results as solutions of both problems 1 and 2. In the case of Problem 2 we
use two values of $k=60,90.$ In particular, we want to evaluate the
resolution of our technique. We define the resolution as the distance
between surfaces of two inclusions at which they can be separated in our
images. Even though distances below are dimensionless, dimensions can be
easily assigned: the dimensionless distance $X$ means $X$ microns.
Thus, our cases are:
\begin{enumerate}
\item Two inclusions of the same size radii, which are symmetric about the $x_2 x_3$ plane: see Figures \ref{fig: two inclusions, distance = 0.30}, \ref{fig: two inclusions, distance = 0.10} and \ref{fig: two inclusions, distance = 0.025}.
\item Two inclusions of the same radii (non-symmetric case): see Figure \ref{fig: two inclusions non symmetric, distance = 0.5}.
\item Two inclusions with different radii: see Figure \ref{fig: two inclusions non equal, distance = 0.5}.
\end{enumerate}

In each Figure, the square has the
center at $\left( x_{1},x_{2}\right) =\left( 0,0\right) $ and its side
equals $\sqrt{2}.$ Hence, this square is inscribed in the circle of the
radius $R=1$ with the center at the origin. We display in each figure: (a)
the 2D cross-section in the $x_{1},x_{2}$ plane of the true image, (b) the
Radon transform $\left( \mathcal{R}\beta \right) \left( r,\theta \right) $
of the true function $\beta \left( x_{1},x_{2},0\right),$ (c) The function $%
l_{B}\left( r,\theta ,60\right) ,$ (d) The reconstructed image $\mathcal{R}%
^{-1}l_{B}\left( r,\theta ,60\right) ,$ (e) The function $l_{B}\left(
r,\theta ,90\right) ,$ (f) The reconstructed image $\mathcal{R}%
^{-1}l_{B}\left( r,\theta ,90\right) ,$ (g) the reconstructed function $\tau
\left( x,x_{0}\right) -\left\vert x-x_{0}\right\vert $ for $x,x_{0}\in
S_{0}\left( R\right) $ in $\left( r,\theta \right) $ coordinates, (h) The approximate
(Remark \ref{Re 4}) and reconstructed phases $\tau \left( x,x_{0}\right) $ when the
point $x$ runs along the bottom side of the square and the source $%
x_{0}=\left( 1,0,0\right) ,$ (i) The approximate
(Remark \ref{Re 4}) and reconstructed moduli $%
\left\vert u\left( x,x_{0},k\right) \right\vert $ of the total field for the
same $x,x_{0}$ as in (h), (j) The reconstructed image via the solution of
Problem 1 (subsection \ref{subsection: 4.2}).

We have
applied a postprocessing procedure to our images. For
each vertex $V$ we found a disk $D\left( V\right) $ with the center at $V$
of the radius $0.005$. Then we have calculated the average value of the imaged
function $\beta $ over all vertices inside of $\overline{D\left( V\right) }.$
Next, we prescribed that average value as the value of $\beta $ at the
vertex $V$. Next, since by (\ref{1.3}) $\beta \geq 0$, we truncated to zero
all negative values of the resulting function $\beta .$ Our images display
functions $\beta $ obtained after this procedure.

It is worth mentioning the software Armadillo by Sanderson \cite{San}. It involves a linear algebra package which is very helpful to speed up our computations.   

\begin{figure}[h!]
		\hfill\subfloat[\labela]{
		\includegraphics[width=\width, height = \height ]{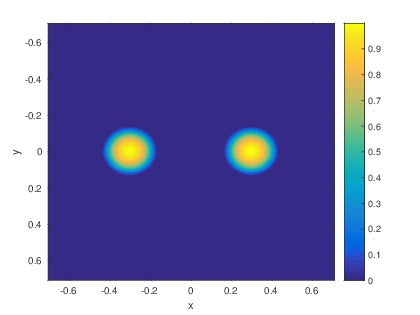} 
	}
	\hfill\subfloat[\labelb]{
		\includegraphics[width=\width, height = \height,height = \height ]{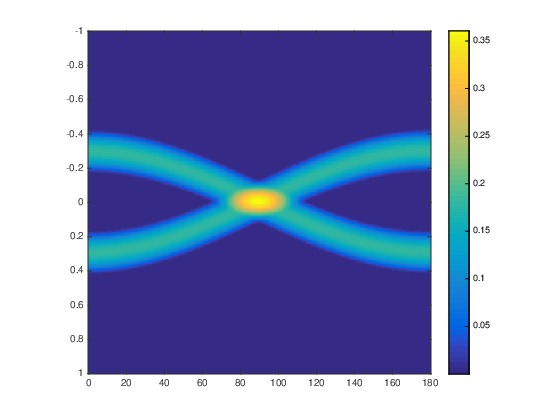} 
	}
	\hfill\subfloat[\labelc]{
	\includegraphics[width=\width, height = \height,height = \height ]{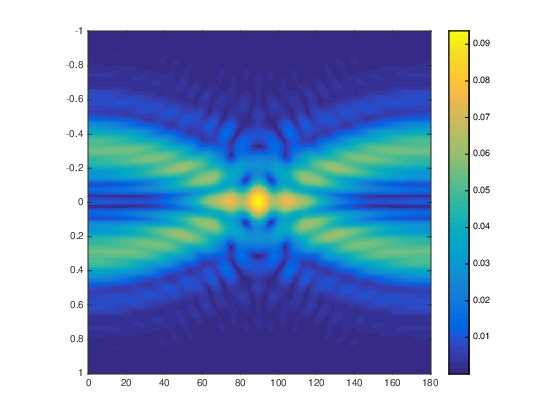} 
	}
	\hfill\subfloat[\labeld]{
	\includegraphics[width=\width, height = \height,height = \height ]{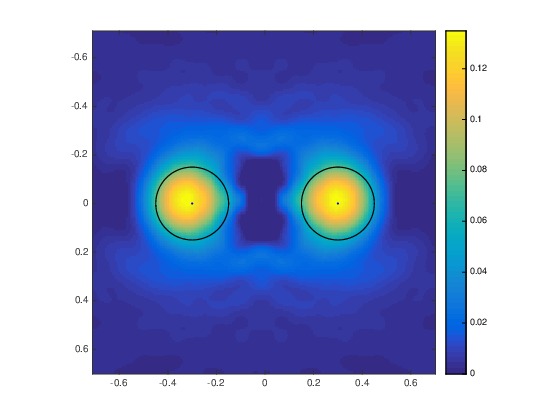} 
	}
	\hfill\subfloat[\labele]{
	\includegraphics[width=\width, height = \height,height = \height ]{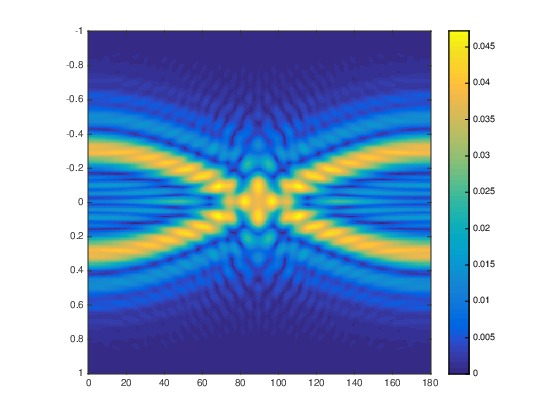} 
	}	  
	\hfill\subfloat[\labelf]{
	\includegraphics[width=\width, height = \height,height = \height]{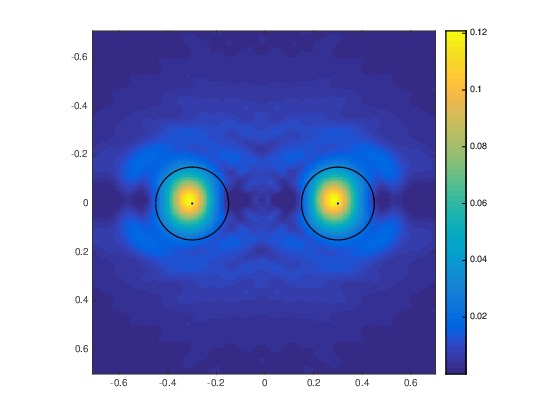} 
	}
	\hfill\subfloat[\labelg]{
		\includegraphics[width=\width, height = \height,height = \height ]{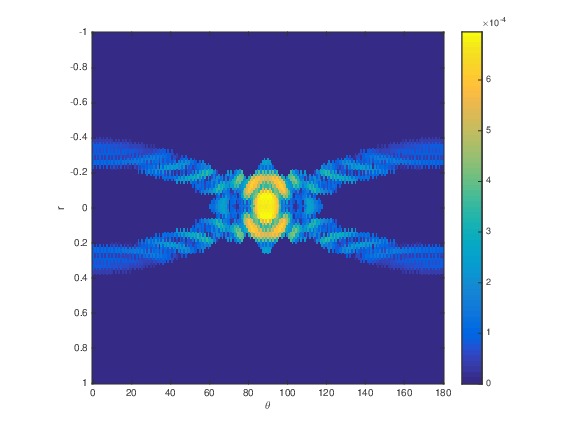} 
	}
	\hfill
	\subfloat[\labelh]{
		\includegraphics[width=\width, height = \height ]{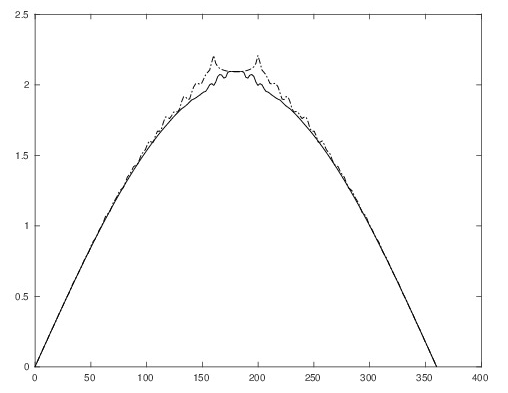} 	
	}
	\hfill
	\subfloat[\labeli]{
		\includegraphics[width=\width, height = \height]{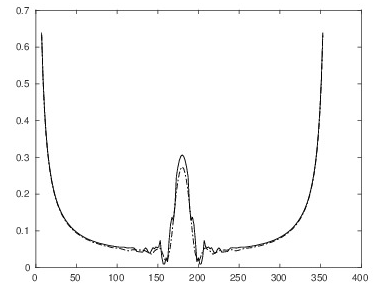} 	
	}
	\hfill\subfloat[\labelk]{
		\includegraphics[width=\width, height = \height,height = \height]{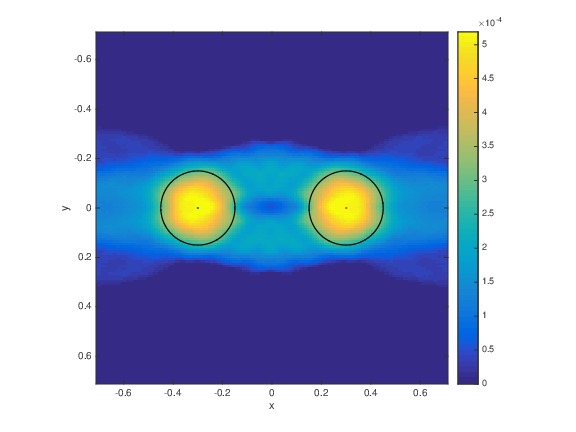} 
	}
	\caption{\label{fig: two inclusions, distance = 0.30} Two inclusions of the radius 0.15 with the distance 0.3 between their
surfaces.}
\end{figure}

\begin{figure}[h!]
		\hfill\subfloat[\labela]{
		\includegraphics[width=\width, height = \height]{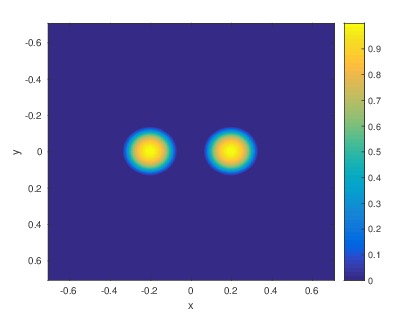} 
	}
	\hfill\subfloat[\labelb]{
		\includegraphics[width=\width, height = \height]{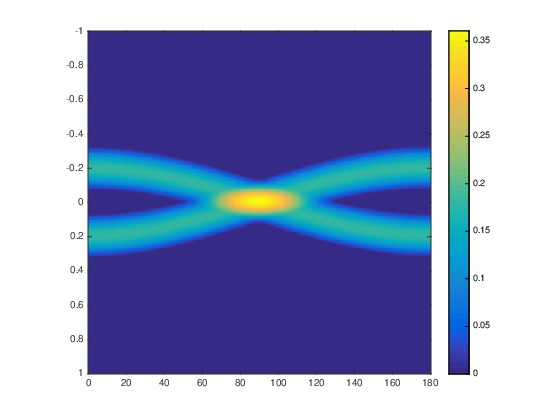} 
	}
	\hfill\subfloat[\labelc]{
	\includegraphics[width=\width, height = \height]{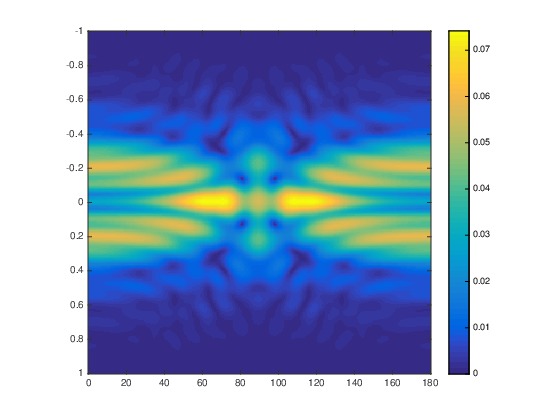} 
	}
	\hfill\subfloat[\labeld]{
	\includegraphics[width=\width, height = \height]{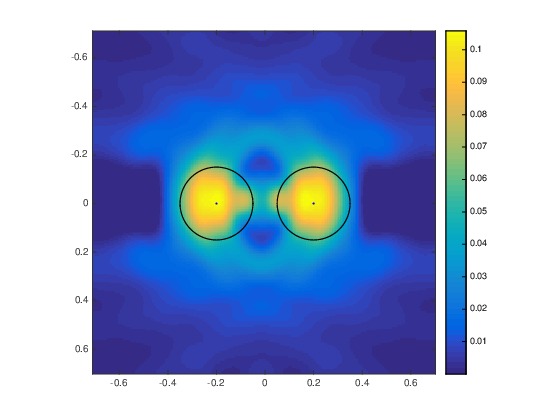} 
	}
	\hfill\subfloat[\labele]{
	\includegraphics[width=\width, height = \height]{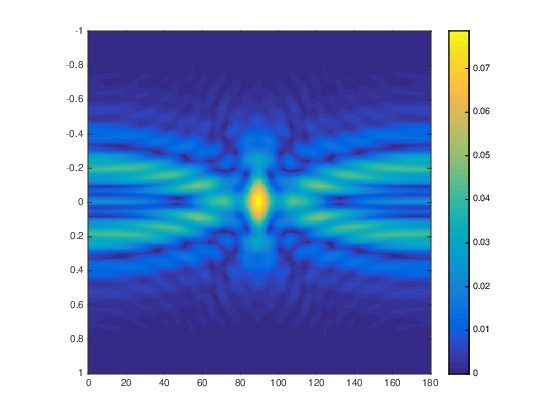} 
	}	  
	\hfill\subfloat[\labelf]{
	\includegraphics[width=\width, height = \height]{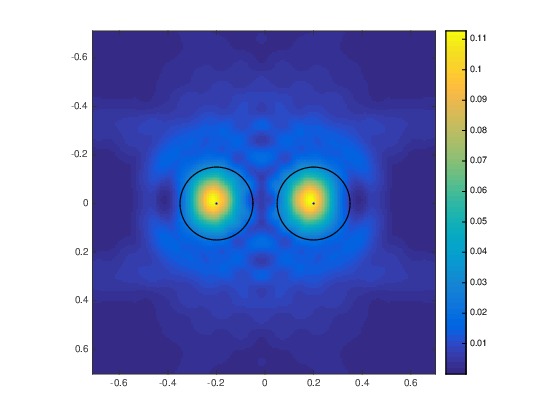} 
	}
	\hfill\subfloat[\labelg]{
		\includegraphics[width=\width, height = \height]{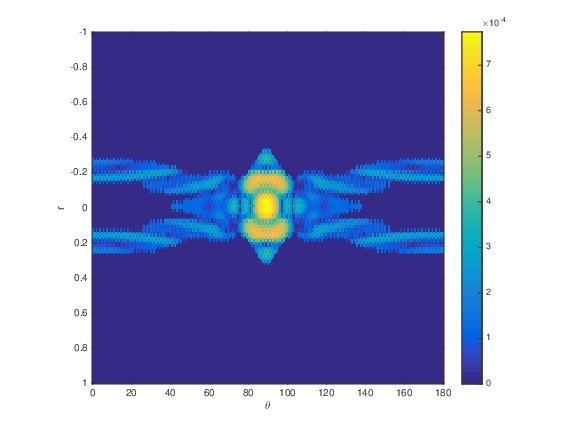} 
	}
	\hfill
	\subfloat[\labelh]{
		\includegraphics[width=\width, height = \height]{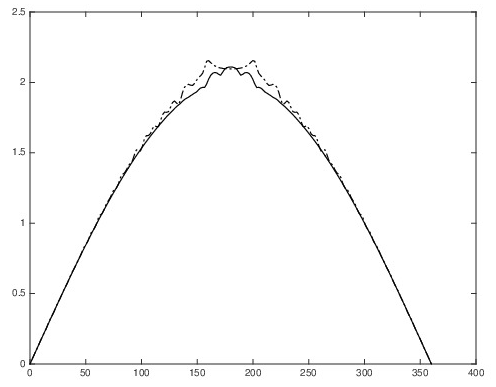} 	
	}
	\hfill
	\subfloat[\labeli]{
		\includegraphics[width=\width, height = \height]{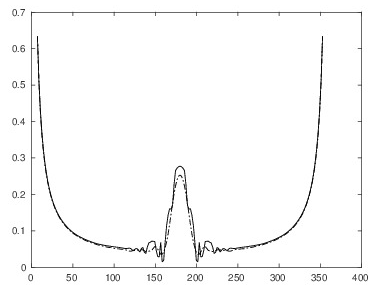} 	
	}
	\hfill\subfloat[\labelk]{
		\includegraphics[width=\width, height = \height]{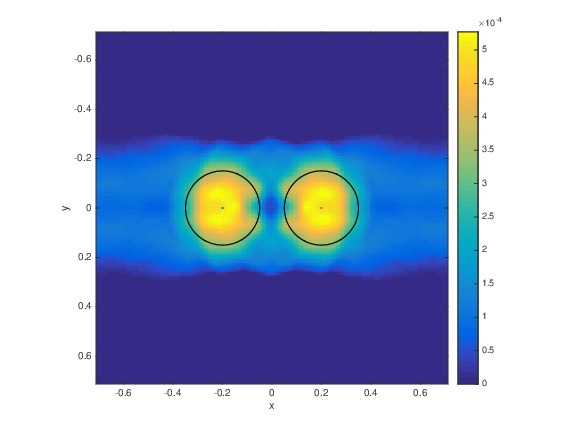} 
	}
	\caption{\label{fig: two inclusions, distance = 0.10} Two inclusions of
the radius 0.15 with the distance 0.1 between their surfaces.}
\end{figure}

\begin{figure}[h!]
	\hfill\subfloat[\labela]{
		\includegraphics[width=\width, height = \height]{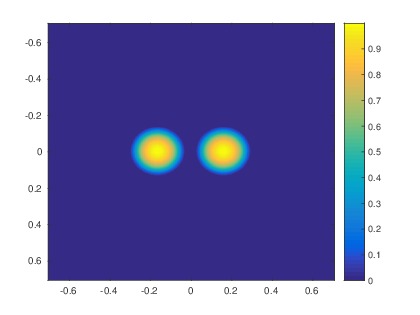} 
	}
	\hfill\subfloat[\labelb]{
		\includegraphics[width=\width, height = \height]{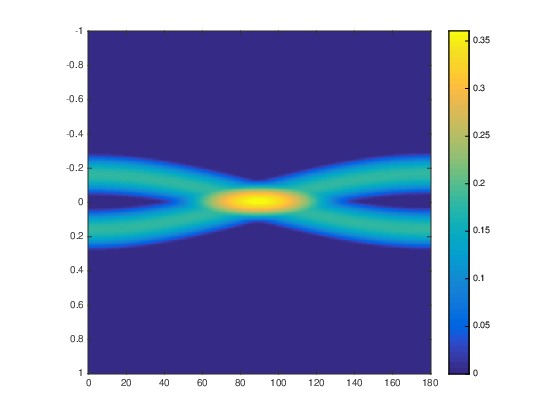} 
	}
	\hfill\subfloat[\labelc]{
	\includegraphics[width=\width, height = \height]{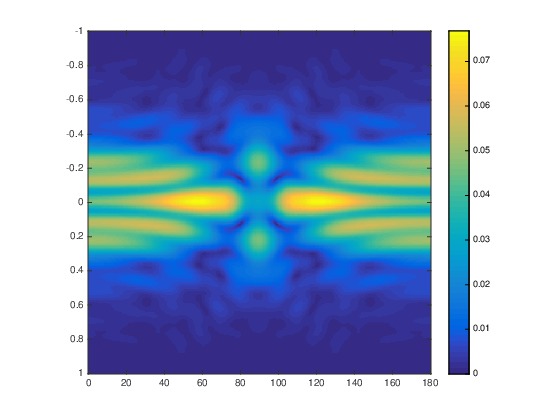} 
	}
	\hfill\subfloat[\labeld]{
	\includegraphics[width=\width, height = \height]{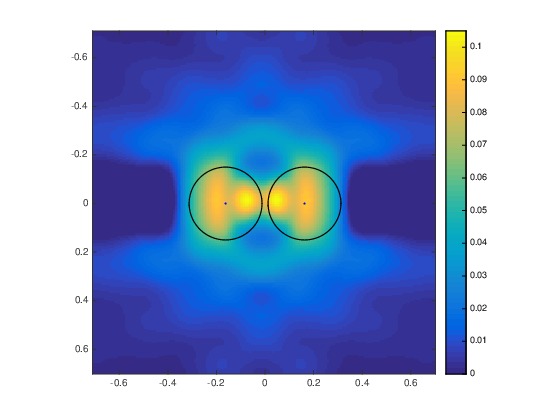} 
	}
	\hfill\subfloat[\labele]{
	\includegraphics[width=\width, height = \height]{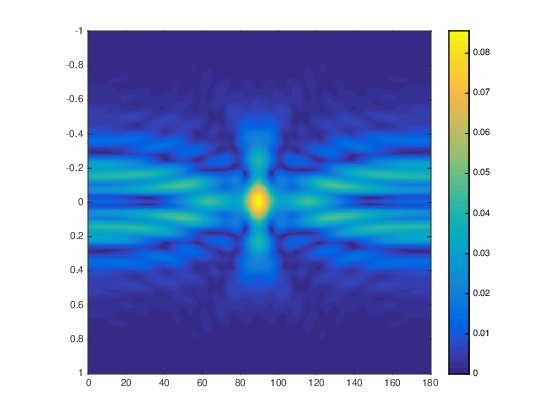} 
	}	  
	\hfill\subfloat[\labelf]{
	\includegraphics[width=\width, height = \height]{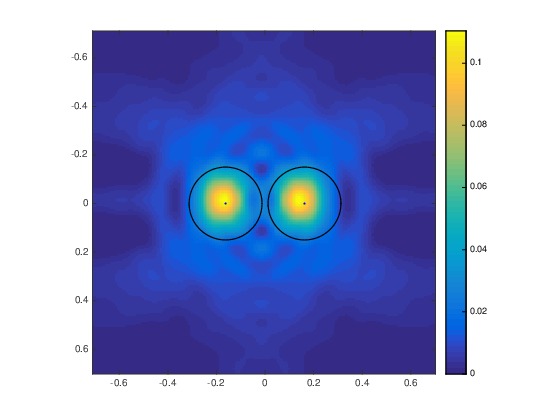} 
	}
	\hfill\subfloat[\labelg]{
		\includegraphics[width=\width, height = \height]{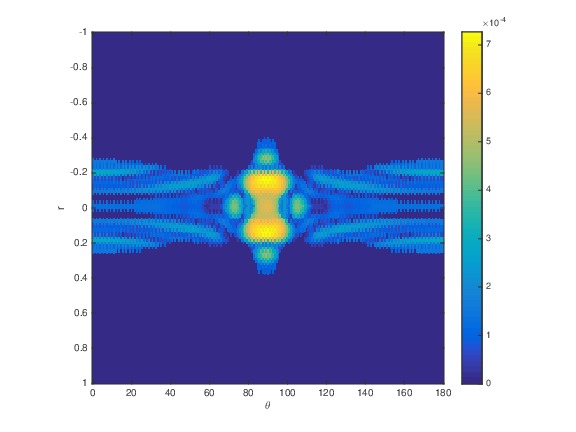} 
	}
	\hfill
	\subfloat[\labelh]{
		\includegraphics[width=\width, height = \height]{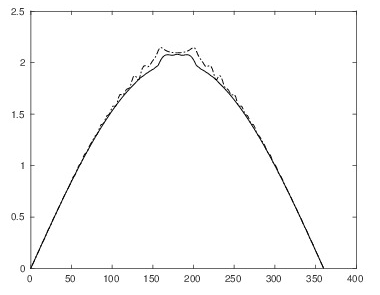} 	
	}
	\hfill
	\subfloat[\labeli]{
		\includegraphics[width=\width, height = \height]{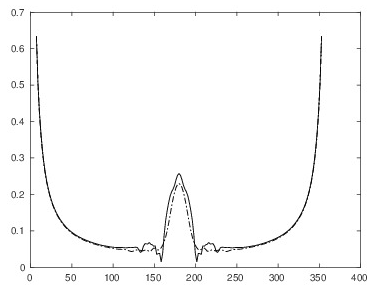} 	
	}
	\hfill\subfloat[\labelk]{
		\includegraphics[width=\width, height = \height]{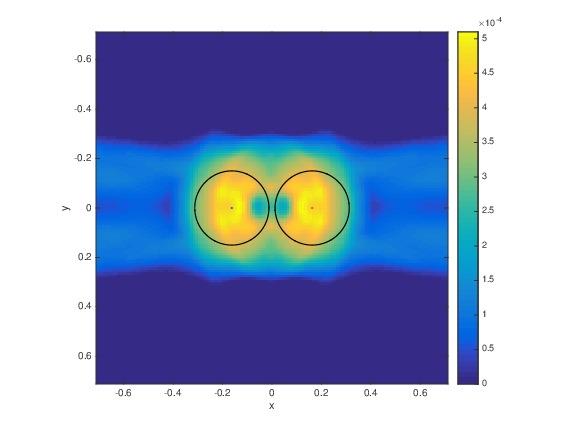} 
	}
	\caption{\label{fig: two inclusions, distance = 0.025} Two inclusions of
the radius 0.15 with the distance 0.025 between their surfaces. This is the
case of super resolution, see section \ref{sec: summ}.}
\end{figure}

\begin{figure}[h!]
		\hfill\subfloat[\labela]{
		\includegraphics[width=\width, height = \height]{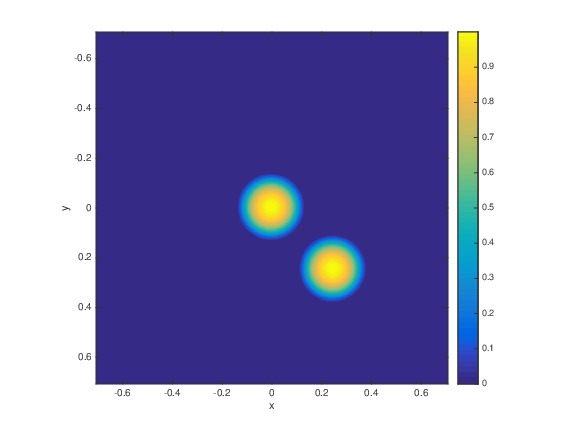} 
	}
	\hfill\subfloat[\labelb]{
		\includegraphics[width=\width, height = \height]{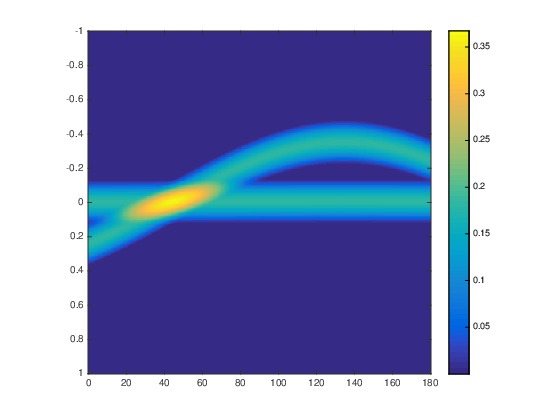} 
	}
	\hfill\subfloat[\labelc]{
	\includegraphics[width=\width, height = \height]{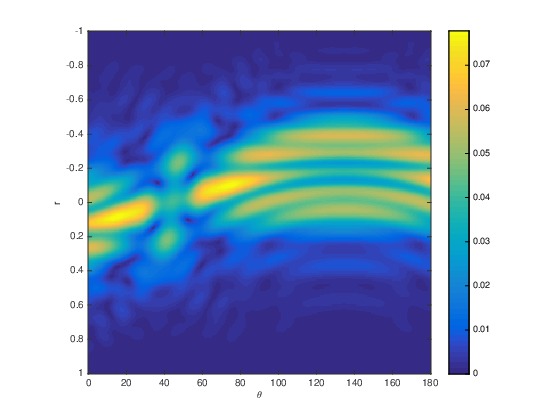} 
	}
	\hfill\subfloat[\labeld]{
	\includegraphics[width=\width, height = \height]{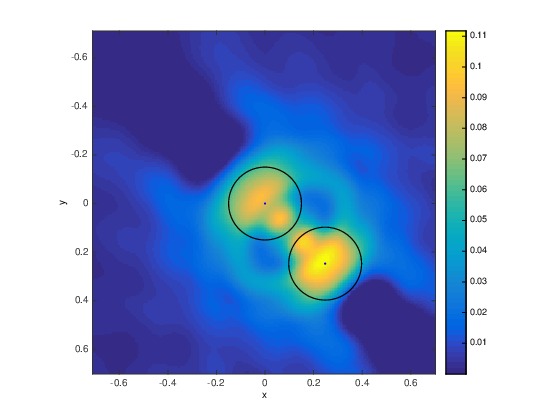} 
	}
	\hfill\subfloat[\labele]{
	\includegraphics[width=\width, height = \height]{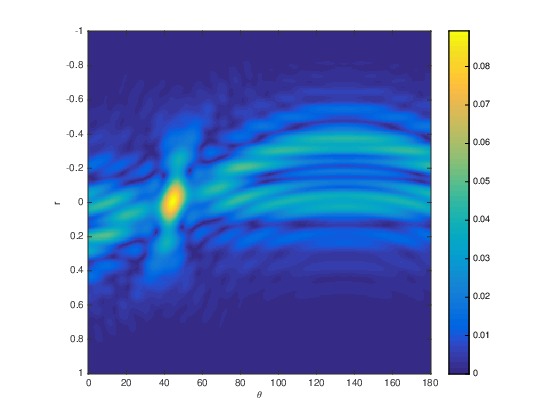} 
	}	  
	\hfill\subfloat[\labelf]{
	\includegraphics[width=\width, height = \height]{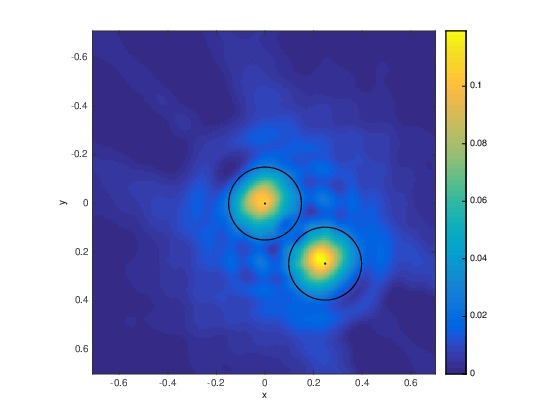} 
	}
	\hfill\subfloat[\labelg]{
		\includegraphics[width=\width, height = \height]{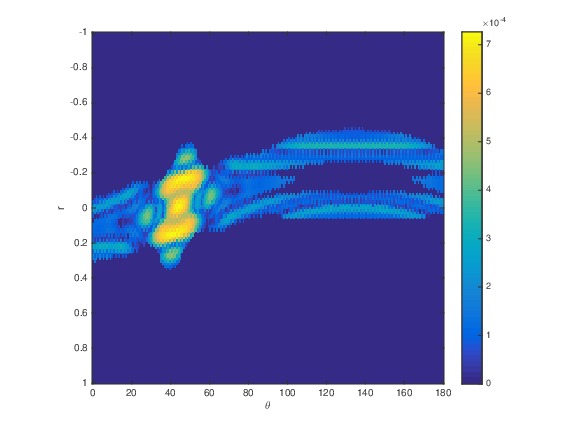} 
	}
	\hfill
	\subfloat[\labelh]{
		\includegraphics[width=\width, height = \height]{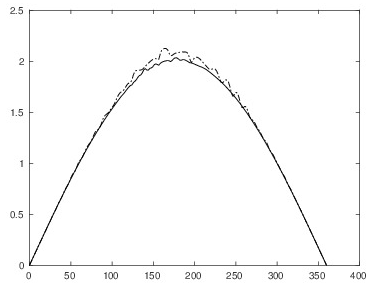} 	
	}
	\hfill
	\subfloat[\labeli]{
		\includegraphics[width=\width, height = \height]{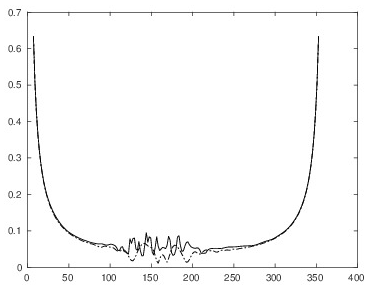} 	
	}
	\hfill\subfloat[\labelk]{
		\includegraphics[width=\width, height = \height]{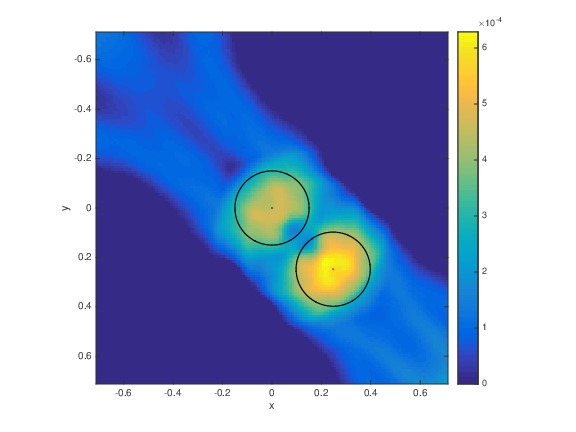} 
	}
	\caption{\label{fig: two inclusions non symmetric, distance = 0.5} Two inclusions of the radius $0.15$ centered at $(0, 0, 0)$ and $(0.2475, 0.2475, 0)$. The distance between their surfaces is 0.05.} 
\end{figure}

\begin{figure}[h!]
	\hfill\subfloat[\labela]{
		\includegraphics[width=\width, height = \height]{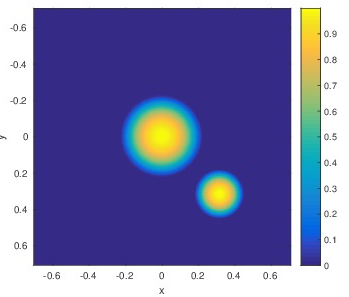} 
	}
	\hfill\subfloat[\labelb]{
		\includegraphics[width=\width, height = \height]{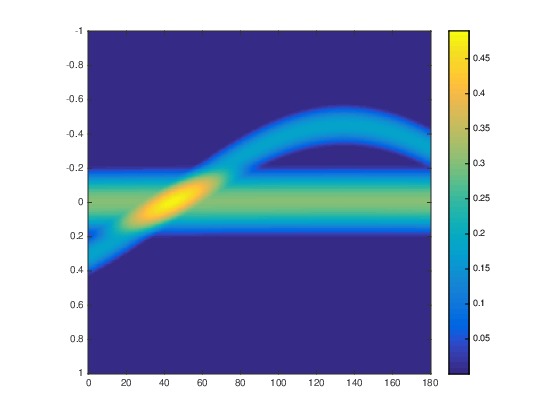} 
	}
	\hfill\subfloat[\labelc]{
	\includegraphics[width=\width, height = \height]{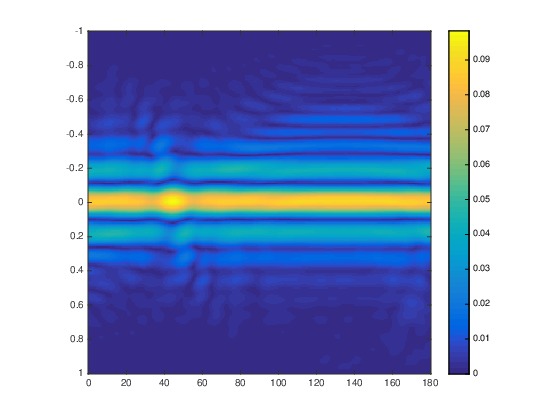} 
	}
	\hfill\subfloat[\labeld]{
	\includegraphics[width=\width, height = \height]{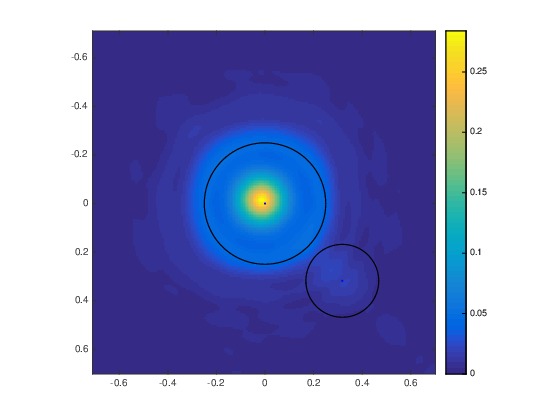} 
	}
	\hfill\subfloat[\labele]{
	\includegraphics[width=\width, height = \height]{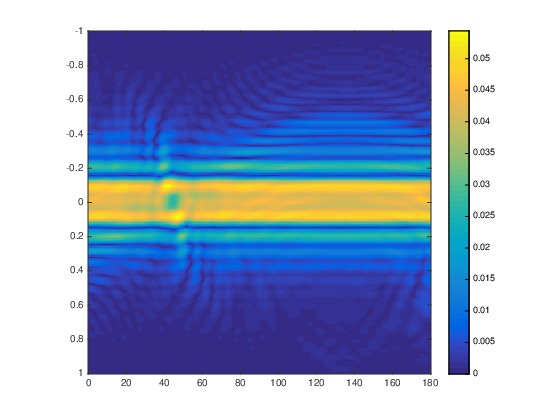} 
	}	  
	\hfill\subfloat[\labelf]{
	\includegraphics[width=\width, height = \height]{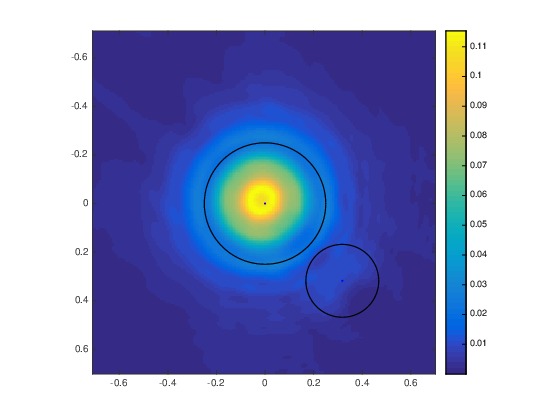} 
	}
	\hfill\subfloat[\labelg]{
		\includegraphics[width=\width, height = \height]{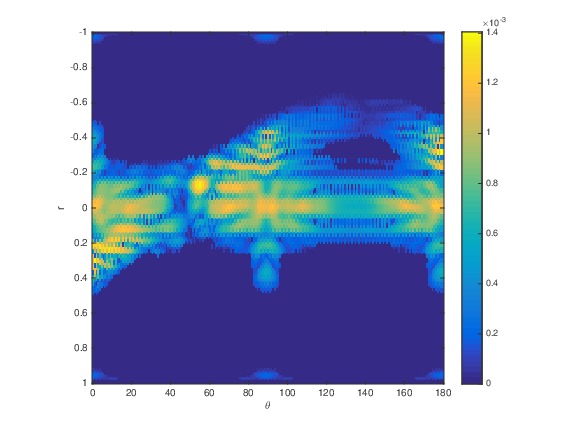} 
	}
	\hfill
	\subfloat[\labelh]{
		\includegraphics[width=\width, height = \height]{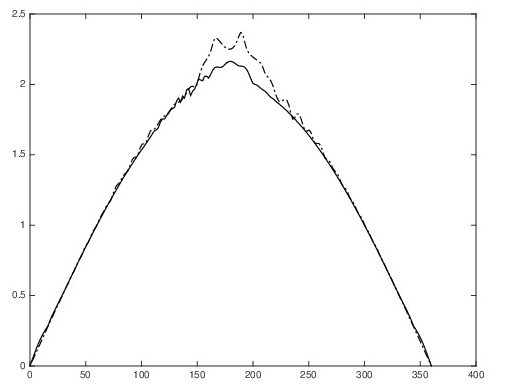} 	
	}
	\hfill
	\subfloat[\labeli]{
		\includegraphics[width=\width, height = \height]{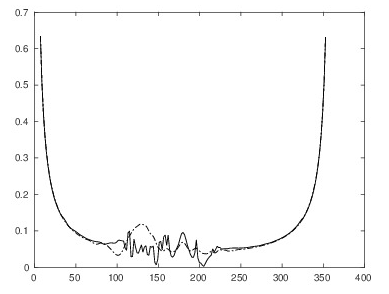} 	
	}
	\hfill\subfloat[\labelk]{
		\includegraphics[width=\width, height = \height]{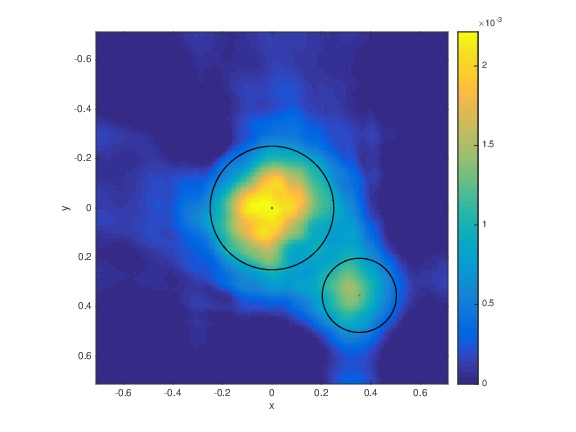} 
	}
	\caption{\label{fig: two inclusions non equal, distance = 0.5} Two inclusions of
the radii 0.25 and 0.15 with the distance 0.05 between their surfaces. The
image of Problem 1 has a better quality than the one of Problem 2.}
\end{figure}

\section{Summary} \label{sec: summ}

In this paper, two rigorous numerical reconstruction procedures are applied
for the first time to solve the 3D inverse scattering problem without the
phase information for the generalized Helmholtz equation (\ref{1.4}). We have conducted our computations
for ranges of parameters, which are realistic
for applications to imaging of nanostructures of hundreds of nanometers
size. Both X-rays and optical signals on the wavelengths range (\ref{101})
can be used as sources of radiation. Imaging of living biological cells
using optical signals is also feasible. 

We have used two reconstruction procedures.\ The first one is based on the idea of 
\cite{KR2}. Here, the specific way of finding the travel time $\tau
\left( x,x_{0}\right) $ for $x,x_{0}\in S_0(R) $ is
significantly different from the one proposed in \cite{KR2}. This is because it
is assumed in \cite{KR2} that a large
interval of frequencies is available. However, by (\ref{1.54}) our realistic dimensionless frequencies deviate, as the maximum, by 30\% ($\left( 65/50-1\right)
\cdot 100\%$) from the central frequency $k_{central}=65$.
The second reconstruction procedure is based on the Born approximation
assumption of \cite{KRB} and it is simpler than the one of \cite{KR2}. 

In both procedures the linearization takes place. While the Born
approximation assumes the linearization from the start, the procedure of 
\cite{KR2} uses the linearization only on the last step: after the function $%
\tau \left( x,x_{0}\right) $ for $x,x_{0}\in S_0(R) $ is found.
The linearization in this case is caused by the fact that it is yet unclear
how to compute numerically the
function $\beta ,$ see, e.g. books \cite{LRV, R1, R2} if the travel time $\tau \left(
x,x_{0}\right) $ is known for $x,x_{0}\in S_0(R).$

Each of our two methods actually assumes that the frequency $k\rightarrow
\infty $ and ends up with the 2D Radon transform. Even though one should
work, in principle at least, with infinite values of $k$, we demonstrate
good quality images for our realistic values of $k$. Furthermore,
even though the Born approximation becomes invalid for large values of $k$,
our computations show that images, using the Born approximation, are
slightly better for $k=90$ than for $k=60$. This is regardless on the fact
that in both cases the data for the inverse problem were generated via the
solution of the corresponding forward problem \eqref{LS eqn} without the Born approximation
assumption. The qualities of images are about the same for both
methods. Still, in the case of Figure \ref{fig: two inclusions, distance = 0.025} the image provided by the solution of Problem
1 has a better quality than the one provided by the Born approximation.

By (\ref{101}) the minimal wavelength we use is $\lambda _{\min }=0.078\mu m.$
On the other hand, we observe on Figure \ref{fig: two inclusions, distance = 0.025} that two inclusions with the
distance of $0.025\mu m=0.32\lambda _{\min }$ are resolved. Since $%
0.32\lambda _{\min }<\lambda _{\min }/2,$ then we observe a rare case of the 
\emph{super resolution}. Another interesting observation for all our figures
\ref{fig: two inclusions, distance = 0.30}-\ref{fig: two inclusions non equal, distance = 0.5} is that our data are quite far from the range of the operator of the
Radon transform: on each of Figures \ref{fig: two inclusions, distance = 0.30}-\ref{fig: two inclusions non equal, distance = 0.5}   compare (b) with (e) and (g).
Nevertheless, the application of the inverse Radon transform to our data
leads to rather good quality images. Although we have not added a random
noise to our data, this difference in the data can also be considered as a
noise in the data.

There is a significant advantage of the method of Problem
1 over the method of Problem 2. Indeed, the method of Problem 1 finds the phase $\tau
\left( x,x_{0}\right) $ as well as the
modulus $A\left( x,x_{0}\right) $ of the wave field $u\left(
x,x_{0},k\right) $ for $x,x_{0}\in S_0(R) $ and for large $k$
with a good accuracy and without any linearization, see (h) and (i) in each of Figures \ref{fig: two inclusions, distance = 0.30}-\ref{fig: two inclusions non equal, distance = 0.5} as well as \eqref{1.110}.
This should pave the way to a future significant image refinement. Indeed,
while we image shapes of abnormalities with a good accuracy, the accuracies of
the calculated abnormalities/background contrasts $\left( 1+\beta \left(
x\right) \right) /1$ are poor. We attribute the latter to the linearization.
On the other hand, as soon as the function $u\left( x,x_{0},k\right) $ is
found accurately for $x,x_{0}\in S_0(R) ,$ one can apply a
modified globally convergent method of \cite{BK} for a corresponding
Coefficient Inverse Problem. The modification is due to the fact that in 
\cite{BK} the time domain data were considered and then the Laplace
transform was applied. Whereas here we would need to have Fourier transform.
An important point to make here is that the method of \cite{BK} images
inclusions/background contrasts accurately. 

Hence, one can apply a two-stage numerical procedure. On the first stage one
would reconstruct the function $u\left( x,x_{0},k\right) $ for $x,x_{0}\in
S_0(R) $ and for large values of $k$. On this stage shapes of
abnormalities will be well reconstructed. On the second stage one would
reconstruct the function $\beta \left( x\right) $ within abnormalities,
while keeping the shape the same as it was obtained on the first stage. On
the second stage, one would use the already calculated approximation for the
function $\beta \left( x\right) $ to find the first approximation for the
so-called \textquotedblleft tail function". This function is an important
element of the method of \cite{BK}.

We note that this is a sort of a reversed two-stage procedure of Chapters 4
and 5 of \cite{BK}.\ In \cite{BK} the globally convergent method is applied
on the first stage to find locations of abnormalities and
inclusions/background contrasts. On the second stage, a sequence of Tikhonov
functionals is minimized in the so-called
\textquotedblleft adaptivity technique" to find the
shapes of inclusions. In doing so, the solution of the first stage is used
as the starting point of iterations. It is important that the method of \cite%
{BK} is completely verified on experimental data, see Chapter 5 of \cite{BK}
 as well as \cite{BK1,BK2,TBKF1,TBKF2}. The
second option for the globally convergent method is the one of \cite{KT},
where a globally convex cost functional is constructed using a Carleman
Weight Function. In both these globally convergent methods only one source
position at a time is used. This means in our case that we can use several
sources sequentially to refine the images further.




\end{document}